\documentclass[11pt,bezier,amstex]{article}  

\usepackage{color}              
\usepackage{graphicx}
\usepackage[]{amsmath}
\usepackage{amssymb}
\usepackage{hyperref}
\usepackage{enumerate, amsfonts, amsthm, bbm}

\definecolor{MyDarkBlue}{rgb}{0,0.08,0.50}
\definecolor{BrickRed}{rgb}{0.65,0.08,0}

\hypersetup{
colorlinks=true,       
    linkcolor=MyDarkBlue,          
    citecolor=BrickRed,        
    filecolor=red,      
    urlcolor=cyan           
}

\newtheorem{Lemma}{Lemma}[section]

\newtheorem{Theorem}[Lemma]{Theorem}

\newtheorem{Corollary}[Lemma]{Corollary}

\newcommand{\prob}{\mathbb{P}}
\newcommand{\bfw}{\boldsymbol{w}}

\newcommand{\CC}{\mathcal{C}}
\newcommand{\FF}{\mathcal{F}}

\newcommand{\GG}{\mathcal{G}}

\newcommand{\eps}{\varepsilon}

\newcommand{\set}[1]{\left\{#1\right\}}

\newcommand{\Rbold}{{\mathbb{R}}}

\newcommand{\ind}[2]{1_{(e \in \pi(#1,#2))}}

\newcommand{\expec}{\mathbb{E}}

\def\ind{{\rm 1\hspace{-0.90ex}1}}

\newcommand{\eqn}[1]{\begin{equation} #1 \end{equation}}
\newcommand{\eqan}[1]{\begin{align} #1 \end{align}}
\newcommand{\lbeq}[1]{\label{#1}}
\newcommand{\refeq}[1]{(\ref{#1})}
\newcommand{\sss}{\scriptscriptstyle}

\newcommand{\op}{o_{\sss \prob}}
\newcommand{\Op}{O_{\sss \prob}}

\newcommand{\Poi}{{\rm Poi}}
\newcommand {\vep}{\varepsilon}
\newcommand\1{\mathbbm{1}}
\newcommand{\indic}[1]{\1_{\{#1\}}}

\newcommand {\convd}{\stackrel{d}{\longrightarrow}}
\newcommand {\convp}{\stackrel{\sss {\mathbb P}}{\longrightarrow}}

\newcommand{\nn}{\nonumber}

\numberwithin{equation}{section}

\begin{document}

\author{Shankar Bhamidi
\thanks{
Department of Statistics and Operations Research,
University of North Carolina,
304 Hanes Hall, Chapel Hill, NC27510, United States.
Email: {\tt shankar@math.ubc.ca}}
\and
Remco van der Hofstad
\thanks{Department of Mathematics and
Computer Science, Eindhoven University of Technology, P.O.\ Box 513,
5600 MB Eindhoven, The Netherlands. E-mail: {\tt
rhofstad@win.tue.nl}, {\tt j.s.h.v.leeuwaarden@tue.nl}}
\and
Johan S.H. van Leeuwaarden$^{~\dagger}$
}

\title{Scaling limits for critical inhomogeneous random graphs with finite third moments}

\maketitle

\begin{abstract}
We identify the scaling limits for the sizes of the largest components at criticality
for inhomogeneous random graphs with weights that have finite third moments.
We show that the sizes of the (rescaled) components converge to the excursion
lengths of an inhomogeneous Brownian motion, which extends results of Aldous \cite{Aldo97} for the critical behavior of Erd\H{o}s-R\'enyi random graphs.
We rely heavily on martingale convergence techniques, and concentration properties
of (super)martingales. This paper is part of a programme initiated in \cite{Hofs09a} to study the near-critical behavior
in inhomogeneous random graphs of so-called rank-1.
\end{abstract}

\vspace{0.6in}


\noindent
{\bf MSC2000 subject classification.}
60C05, 05C80, 90B15.

\section{Introduction}
\label{sec-int}

\subsection{Model}
\label{sec-mod}

We start by describing the model considered in this paper.
While there are many variants available in the
literature, the most convenient for our purposes model is the model often referred to as \emph{Poissonian graph process}  or  \emph{Norros-Reittu model} \cite{NorRei06}. See Section \ref{sec-disc} below for consequences
for other models. To define the model, we consider the vertex set
$[n]:=\{1,2,\ldots, n\}$, and attach an edge with probability $p_{ij}$
between vertices $i$ and $j$, where
    \eqn{
    \lbeq{pij-def}
    p_{ij} =1-\exp\left( -\frac{w_i w_j}{l_n}\right),
    }
and
    \eqn{
    \lbeq{ln-def}
    l_n = \sum_{i=1}^n w_i.
    }
Different edges are independent.

Below, we shall formulate general conditions on the weight sequence
$\bfw = (w_1,\ldots, w_n)$, and for now formulate two main examples.
The first key example arises when we take
$\bfw$ to be an i.i.d.\ sequence of random variables
with distribution function $F$ satisfying
    \eqn{
    \lbeq{[1-F]bd}
    \expec[W^3]<\infty,
    }
The second key example (which is also studied in \cite{Hofs09a}) arises when we let the weight sequence $\bfw=(w_i)_{i=1}^n$
be defined by
    \eqn{
    \lbeq{choicewi}
    w_i = [1-F]^{-1}(i/n),
    }
where $F(x)$ is a distribution function satisfying that, with $W$ a random variable
with distribution function $F$ satisfying \refeq{[1-F]bd},
and where $[1-F]^{-1}$ is the generalized inverse function of
$1-F$ defined, for $u\in (0,1)$, by
    \eqn{
    \lbeq{invverd}
    [1-F]^{-1}(u)=\inf \{ s: [1-F](s)\leq u\}.
    }
By convention, we set $[1-F]^{-1}(1)=0$.

Write
    \eqn{
    \lbeq{nu-def}
    \nu = \frac{\expec[W^2]}{\expec[W]}.
    }
Then, by \cite{BolJanRio07}, the random graphs we consider are subcritical when $\nu<1$ and
supercritical when $\nu>1$. Indeed, when $\nu> 1$, then there is one giant component of
size $\Theta(n)$ while all other components are of smaller size $\op(n)$, while when
$\nu \leq 1$ the largest connected component has size $\op(n)$.
Thus, the critical value of the model is $\nu=1$. Here, and throughout this paper, we use the following standard notation.
We write $f(n)=O(g(n))$ for functions $f,g\ge0$ and $n\rightarrow\infty$ if there
exists a constant $C>0$ such that $f(n)\le Cg(n)$ in the limit, and $f(n)=o(g(n))$
if $g(n)\neq O(f(n))$. Furthermore, we write  $f=\Theta(g)$ if $f=O(g)$ and $g=O(f)$.
We write $\Op(b_n)$ for a sequence of random variables $X_n$ for which $|X_n|/b_n$
is tight as $n\rightarrow\infty$, $\op(b_n)$ for a sequence of random variables $X_n$
for which $|X_n|/b_n\convp0$ as $n\rightarrow\infty$.
Finally, we write that a sequence of events $(E_n)_{n\geq 1}$
occurs \emph{with high probability} ({\bf whp}) when $\prob(E_n)\rightarrow 1$.

We shall write $\GG_n^0(\bfw)$ to be the
graph constructed via the above procedure, while, for any fixed $t\in \Rbold$, we
shall write $\GG_n^t(\bfw)$ when we use the weight sequence $(1+t n^{-1/3})\bfw$, for which the
probability that $i$ and $j$ are neighbors equals $1-\exp\left(-(1+t n^{-1/3})w_iw_j/l_n\right)$. In this setting we take $n$ so large that $1+t n^{-1/3}>0$.

We now formulate the general conditions on the weight sequence
$\bfw$. In Section \ref{sec-ver-cond},
we shall verify that these conditions are satisfied for i.i.d.\ weights with
finite third moment, as well as for the choice in \refeq{choicewi}.
We assume the following three conditions on the weight sequence $\bfw$:

\paragraph{(a) Maximal weight bound.}
We assume that the maximal weight is
$o(n^{1/3})$, i.e.,
    \eqn{
    \lbeq{max-bd-cond}
    \max_{i\in [n]} w_i=o(n^{1/3}).
    }
\paragraph{(b) Weak convergence of weight distribution.} We assume that
the weight of a uniform vertex converges in distribution to some
distribution function $F$, i.e., let $V_n\in [n]$ be a uniform vertex. Then we assume that
    \eqn{
    \lbeq{weak-conv-cond}
    w_{V_n}\convd W,
    }
for some limiting random variable $W$ with distribution function $F$.
Condition \eqref{weak-conv-cond} is equivalent to the statement that,
for every $x$ which is a continuity point of $x\mapsto F(x)$, we have
    \eqn{
    \lbeq{weak-conv-cond-eq}
    \frac{1}{n}\#\{i: w_i\leq x\}\rightarrow F(x).
    }
\paragraph{(c) Convergence of first three moments.} We assume that
    \eqan{
    \frac{1}{n}\sum_{i=1}^n w_i&=\expec[W]+o(n^{-1/3}),\lbeq{conv-first-mom}\\
    \frac{1}{n}\sum_{i=1}^n w_i^2&=\expec[W^2]+o(n^{-1/3}),\lbeq{conv-sec-mom}\\
    \frac{1}{n}\sum_{i=1}^n w_i^3&=\expec[W^3]+o(1)\lbeq{conv-third-mom}.
    }
Note that condition (a) follows from conditions (b) and (c), as
we prove around \refeq{unif-int-w3} below.
When $\bfw$ is \emph{random}, for example in the case where $(w_i)_{i=1}^n$
are i.i.d.\ random variables with finite third moment, then we need the
estimates in conditions (a), (b) and (c) to hold \emph{in probability}.

We shall simply refer to the above three conditions as
conditions (a), (b) and (c). Note that \refeq{conv-first-mom} and
\refeq{conv-sec-mom} in condition (c) also imply that
    \eqn{
    \lbeq{nu_n-asymp}
    \nu_n=\frac{\sum_{i=1}^n w_i^2}{\sum_{i=1}^n w_i}=\frac{\expec[W^2]}{\expec[W]}+o(n^{-1/3})
    =\nu+o(n^{-1/3}).
    }

Before we write our main result we shall need one more construct.
For fixed $t\in \Rbold$ consider the inhomogeneous Brownian motion $(W^t(s))_{s\geq 0}$ with
    \begin{equation}
	W^t(s) =  B(s) + st -\frac{s^2}{2},
	\label{eqn:wtaldous}
    \end{equation}
where $B$ is standard Brownian motion, and that has drift $t-s$ at time $s$.
We want to consider this process
restricted to be non-negative, which is why we introduce the reflected process
    \begin{equation}
	\bar{W}^t(s) = W^t(s) - \min_{0\leq s^\prime \leq s} W^t(s^\prime).
	\label{eqn:bt}
	\end{equation}
It \cite{Aldo97} it is shown that the excursions of $\bar{W}^t$ from $0$ can
be ranked in increasing order as, say, $\gamma_1(t)> \gamma_2(t)> \ldots$.

Now let $\CC_n^1(t)\geq \CC_n^2(t) \geq \CC_n^3(t)\ldots $ denote the sizes
of the components in $\GG_n^t(\bfw)$ arranged in increasing order.
Define $l^2$ to be the set of infinite sequences
$x=(x_i)_{i=1}^{\infty}$ with $x_1\geq x_2\geq \ldots\geq 0$ and
$\sum_{i=1}^{\infty} x_i^2<\infty$, and define the
$l^2$ metric by
    \eqn{
    d(x,y) =\Big(\sum_{i=1}^{\infty}(x_i-y_i)^2\Big)^{1/2}.
    }
Let
\eqn{
    \mu=\expec[W], \quad  \sigma_3=\expec[W^3].
    }
Then, our main result is as follows:
\begin{Theorem}[The critical behavior]
\label{theo:main}
Assume that the weight sequence $\bfw$
satisfies conditions {\rm(a)}, {\rm(b)} and {\rm(c)}. Then, as $n\to\infty$,
    \eqn{
    \label{resultoex}
    \big(n^{-2/3}\CC_n^i(t)\big)_{i\geq 1} \convd \big(\mu\sigma_3^{-1/3}\gamma_i(t\mu \sigma_3^{-2/3})\big)_{i\geq 1}=:\big(\gamma_i^*(t)\big)_{i\geq 1},
    }
in distribution and with respect to the $l^2$ topology.
\end{Theorem}

Theorem \ref{theo:main} extends the work of Aldous \cite{Aldo97}, who identified the scaling limit of
the largest connected components in the Erd\H{o}s-R\'enyi random graph. Indeed, he proved for the critical Erd\H{o}s-R\'enyi random graph with $p=(1+tn^{-1/3})/n$ that the ordered connected components are given by $\left(\gamma_i(t)\right)_{i\geq 1}$, i.e., the ordered excursions of the reflected process in \eqref{eqn:bt}. Hence, Aldous' result corresponds to Theorem \ref{theo:main}
with $\mu=\sigma_3=1$. The sequence $\big(\gamma_i^*(t)\big)_{i\geq 1}$ is in fact the sequence of ordered excursions of the reflected version of the process
 \begin{equation}
	W^t_*(s) = \sqrt{\frac{\sigma_3}{\mu}} B(s) + st -\frac{s^2\sigma_3}{2\mu^2},
	\label{eqn:wt}
\end{equation}
which reduces to the process in \eqref{eqn:wtaldous} again when $\mu=\sigma_3=1$.

We next investigate the two key examples, and show that conditions (a), (b) and (c)
indeed hold in this case:

\begin{Corollary}[Theorem \ref{theo:main} holds for key examples]
\label{cor:main}
Conditions {\rm(a)}, {\rm(b)} and {\rm(c)} are satisfied in the case where
$\bfw$ is either as in {\rm \refeq{choicewi}}
where $F$ is a distribution function of a random variable $W$ with
$\expec[W^3]<\infty$, or when $\bfw$ consists of
i.i.d.\ copies of a random variable $W$ with $\expec[W^3]<\infty$.
\end{Corollary}

Theorem \ref{theo:main}, in conjunction with
Corollary \ref{cor:main}, proves \cite[Conjecture 1.6]{Hofs09a}, where
the result in Theorem \ref{theo:main} is conjectured in the case where
$\bfw $ is as in \refeq{choicewi}, where
$F$ is a distribution function of a random variable $W$ with
$\expec[W^{3+\vep}]<\infty$ for some $\vep>0$. The current result
implies that $\expec[W^3]$ is a sufficient condition
for this result to hold, and we believe this condition
also to be necessary (as the constant $\expec[W^3]$ also
appears in our results, see \eqref{resultoex} and \eqref{eqn:wt}). Note, however, that
in \cite[Conjecture 1.6]{Hofs09a}, the constant in front of $-s^2/2$ in
\eqref{eqn:wt} is erroneously taken as 1, while it should be
$\expec[W^3]/\expec[W]^2$.

\subsection{Overview of the proof}
\label{sec-cluster-exploration}
In this section, we give an overview of the proof of Theorem \ref{theo:main}.
After having set the stage for the proof, we shall provide a heuristic
that indicates how our main result comes about. We start by describing the
cluster exploration:

\paragraph{Cluster exploration.}
The proof involves two key ingredients:
\begin{itemize}
	\item The exploration of components via breath-first search; and
	\item The labeling of vertices in a size-biased order of their weights $\bfw$.
\end{itemize}
More precisely, we shall explore components and simultaneously construct
the graph $\GG_n^t(\bfw)$ in the following manner: First, for all
ordered pairs of vertices $(i,j)$, let $V(i,j)$ be exponential random variables with rate $\left(1+tn^{-1/3}\right)w_j/l_n$ random variables.
Choose vertex $v(1)$ with probability proportional to $\bfw$, so that
    \eqn{
    \prob(v(1)=i)= w_i/l_n.
    }
The children of $v(1)$ are all the vertices $j$ such that
    \eqn{
    V(v(1),j) \leq w_{v(1)}.
    }
Suppose $v(1)$ has $c(1)$ children. Label these as $v(2), v(3), \ldots v(c(1)+1)$
in increasing order of their $V(v(1),\cdot)$ values. Now move on to $v(2)$ and
explore all of its children (say $c(2)$ of them) and label them as before. Note that when we explore the children of $v(2)$, its potential children cannot include the vertices that we have already identified. More precisely, the children of $v(2)$ consists of the set
    \[\{v\notin \{v(1), \ldots v(c(1)+1)\}: V(v(2),v) \leq w_{v(2)}\}\]
and so on.
Once we finish exploring one component, we move onto the next component
by choosing the starting vertex in a size-biased manner amongst remaining
vertices and start exploring its component. It is obvious that this
constructs all the components of our graph $\GG_n^t(\bfw)$.

Write the breadth-first walk associated to this exploration process as
    \eqn{
    \lbeq{Zn-def}
    Z_n(0) = 0, \hspace{4mm} Z_n(i) = Z_n(i-1) + c(i)-1,
    }
for $i=1,\ldots, n$. Suppose $\CC^*(i)$ is the size of the $i{\rm th}$
component explored in this manner (here we write $\CC^*(i)$ to distinguish
this from $\CC_n^i(t)$, the $i{\rm th}$ {\it largest} component). Then
these can be easily recovered from the above walk by the following prescription:
For $j\geq 0$, write $\eta(j)$ as the stopping time
    \eqn{
    \eta(j) = \min\{i: Z_n(i) = -j\}.
    }
Then
    \eqn{
    \CC^*(j) = \eta(j) - \eta(j-1).
    }
Further,
    \begin{equation}
    Z_n(\eta(j)) = -j, \hspace{4mm} Z_n(i) \geq -j \mbox{ for all }
    \eta(j)< i< \eta(j+1).
    \label{eqn:walk-comp}
    \end{equation}
Recall that we started with vertices labeled $1,2,\ldots, n$
with corresponding weights $\bfw =(w_1, w_2, \ldots, w_n)$.
The size-biased order $v^*(1), v^*(2), \ldots, v^*(n)$ is a
random reordering of the above vertex set where $v(1)=i$ with
probability equal to $w_i/l_n$. Then, given $v^*(1)$, we have that
$v^*(2) = j \in [n]\setminus\{v^*(1)\}$ with probability proportional
to $w_j$ and so on. By construction and the properties of the exponential
random variables, we have the following representation, which lies
at the heart of our analysis:
\begin{Lemma}[Size-biased reordering of vertices]
\label{lemconstruct}
The order $v(1), v(2), \ldots, v(n)$ in the above construction of the
breadth-first exploration process is the size-biased ordering  $v^*(1), v^*(2), \ldots, v^*(n)$ of the
vertex set $[n]$ with weights proportional to $\bfw$.
\end{Lemma}

\proof The first vertex $v(1)$ is chosen from $[n]$ via the size-biased distribution. Suppose it has no neighbors. Then, by construction, the next vertex is chosen via the size-biased distribution amongst all remaining vertices. If vertex $1$ does have neighbors, then we shall use the following construction.

For $j\geq 2$ choose $\tau_{1j}$ exponentially distributed with rate $(1+tn^{-1/3})w_j/l_n$. Rearrange the vertices in increasing order of their $\tau_{1j}$ values (so that $v^\prime(2)$ is the vertex with the smallest $\tau_{1j}$ value, $v^\prime(3)$ is the vertex with the second smallest value and so on). Note that by the properties of the exponential distribution
    \eqn{\prob(v(2) = i \mid v(1))
    = \frac{w_i}{\sum_{j\neq v(1)} w_j}
    \hspace{4mm} \mbox{for } j\in [n]\setminus v(1).}
Similarly, given the value of $v(2)$,
    \eqn{
    \prob(v(3) =i \mid v(1), v(2)) = \frac{w_i}{\sum_{j\neq v(1), v(2)} w_j},
    }
and so on. Thus the above gives us a size-biased ordering of the vertex sex $[n]\setminus \{v(1)\}$. Suppose $c(1)$ of the exponential random variables are less than $w_1$. Then set $v(j) = v^\prime(j)$ for $2\leq j\leq c(1)+1$ and discard all the other labels. This gives us the first $c(1)+1$ values of our size-biased ordering.

Once we are done with $v(1)$, let the potentially unexplored neighbors of $v(2)$ be
    \eqn{
    \mathcal{U}_2 = [n] \setminus \{v(1), \ldots v(c(1)+1)\},
    }
and, again, for $j$ in $\mathcal{U}_2$, we let $\tau_{2j}$ be exponential with rate $(1+tn^{-1/3}) w_j/l_n$ and proceed as above.

Proceeding this way, it is clear that at the end, the random ordering
$v(1), v(2), \ldots, v(n)$ that we obtain is a size-biased random ordering
of the vertex set $[n]$. This proves the lemma.
\qed

\paragraph{Heuristic derivation of Theorem \ref{theo:main}.}
We next provide a heuristic that explains the limiting process in
\eqref{eqn:wt}. Note that by our assumptions on the weight sequence, for the graph $\GG_n^t(\bfw)$
    \eqn{p_{ij} = \left(1+o(n^{-1/3})\right)p^*_{ij},\label{eqn:pij-pijst}}
where
    \eqn{p^{*}_{ij} = \left(1+tn^{-1/3}\right) \frac{w_i w_j}{l_n}.
    }
In the remainder of the proof, wherever we need $p_{ij}$, we shall use $p_{ij}^*$ instead, which shall simplify the calculations and exposition.

Recall the cluster exploration described above, and,
in particular, Lemma \ref{lemconstruct}.
We explore the cluster one vertex at a time, in breadth-first
search. We choose $v(1)$ according to $\bfw$, i.e.,
$\prob(v(1)=j)=w_j/l_n$. We say that a vertex is
\emph{explored} when its neighbors have been investigated,
and \emph{unexplored} when it has been found to be part
of the cluster found so far, but its neighbors have not been
investigated yet. Finally, we say that a vertex is \emph{neutral},
when it has not been considered at all.
Thus, in our cluster exploration, as long as there are unexplored vertices,
we explore the vertices $(v(i))_{i\geq 1}$ in the order of appearance.
When there are no unexplored vertices left, then we draw (size-biased)
from the neutral vertices. Then, Lemma \ref{lemconstruct}
states that $(v(i))_{i=1}^n$ is a size-biased reordering
of $[n]$.

Let $c(i)$ denote the number of neutral neighbors of $v(i)$, and
denote the process $(Z_n(l))_{l\in [n]}$ by $Z_n(l)=0$ and $Z_n(l)=Z_n(l-1)+c(l)-1$.
The clusters of our random graph are found in between successive times in which
$(Z_n(l))_{l\in [n]}$ reaches a new minimum.
Now, Theorem \ref{theo:main} follows from the fact that
$\bar{Z}_n(s) = n^{-1/3} Z_n(\lfloor n^{2/3}s\rfloor)$ weakly converges to
$(W^t_*(s))_{s\geq 0}$ defined as in \eqref{eqn:wt}.
General techniques from
\cite{Aldo97} show that this also implies that the ordered excursions between
successive minima of $(\bar{Z}_n(s))_{s\geq 0}$ also converge to the ones
of $(W^t_*(s))_{s\geq 0}$. These ordered excursions were denoted by
$\gamma_1^*(t)>\gamma_2^*(t)>\ldots$.
Using Brownian scaling, it can be seen that
 \eqn{
    \label{scalrelt}
    W^t_*(s)\stackrel{d}{=}
    \sigma_3^{1/3}\mu^{-2/3}W^{t\mu\sigma_3^{-2/3}}(\sigma_3^{1/3}\mu^{-1}s)
    }
with $W_t$ defined in \eqref{eqn:bt}.
Hence, from the relation \eqref{scalrelt} it immediately follows that
\eqn{
    \lbeq{eq-scal}
    \big(\gamma_i^*(t)\big)_{i\geq 1}\stackrel{d}{=}
    \big(\mu\sigma_3^{-1/3}\gamma_i(t\mu \sigma_3^{-2/3})\big)_{i\geq 1},
    }
which then proves Theorem \ref{theo:main}.

To see how to derive \eqref{scalrelt}, fix $a>0$ and note that $(B(a^2s))_{s\geq 0}$
has the same distribution as $(a B(s))_{s\geq 0}$.
Thus, we obtain, for $(W_{\sigma,\kappa}^t(s))_{s\geq 0}$
with
    \eqn{
    W_{\sigma,\kappa}^t(s)=\sigma B(s)+st-\kappa s^2/2,
    }
the scaling relation
    \eqn{
    W_{\sigma,\kappa}^t(s)\stackrel{d}{=}
    \frac{\sigma}{a} W_{1,a^{-3}\kappa/\sigma}^{t/(a\sigma)}(a^2s).
    }
Using $\kappa=\sigma^2/\mu$ and $a=(\kappa/\sigma)^{1/3}=(\sigma/\mu)^{1/3}$, we note that
   \eqn{
   W_{\sigma,\sigma^2/\mu}^t(s)\stackrel{d}{=}
    \sigma^{2/3}\mu^{1/3}W^{t\sigma^{-4/3}\mu^{1/3}}(\sigma^{2/3}\mu^{-2/3}s),
    }
which, with $\sigma=(\sigma_3/\mu)^{1/2}$ yields \eqref{scalrelt}.

We complete the sketch of proof by giving a heuristic argument that indeed
$\bar{Z}_n(s) = n^{-1/3} Z_n(\lfloor n^{2/3}s\rfloor)$ weakly converges to
$(W^t_*(s))_{s\geq 0}$. For this, we investigate $c(i)$, the number of neutral neighbors
of $v(i)$. Throughout this paper, we shall denote $\tilde w_j=w_j(1+tn^{-1/3})$,
so that the $\GG_n^t(\bfw)$ has weights $\tilde \bfw=(\tilde w_j)_{j\in[n]}$.

We note that since $p_{ij}$ in   \eqref{pij-def} is quite small, the number of neighbors
of a vertex $j$ is close to $\Poi(\tilde w_j)$, where $\Poi(\lambda)$ denotes a Poisson
random variable with mean $\lambda$. Thus, the number of neutral neighbors is
close to the total number of neighbors minus the active neighbors, i.e.,
    \eqn{
    c(i)\approx \Poi(\tilde w_{v(i)})-\Poi\Big(\sum_{j=1}^{i-1} \frac{\tilde w_{v(i)}\tilde w_{v(j)}}{l_n}\Big),
    }
since $\sum_{j=1}^{i-1} \frac{\tilde w_{v(i)}\tilde w_{v(j)}}{l_n}$ is,
conditionally on $(v(j))_{j=1}^i$, the expected number of edges between
$v(i)$ and $(v(j)_{j=1}^{i-1}$. We conclude that the increase of the process $Z_n(l)$ equals
    \eqn{
    c(i)-1\approx \Poi(\tilde w_{v(i)})-1-\Poi\Big(\sum_{j=1}^{i-1} \frac{\tilde w_{v(i)}\tilde w_{v(j)}}{l_n}\Big),
    }
so that
    \eqn{
    \lbeq{Zl-asymp}
    Z_n(l)\approx\sum_{i=1}^l (\Poi(\tilde w_{v(i)})-1)-\Poi\Big(\sum_{j=1}^{i-1} \frac{\tilde w_{v(i)}\tilde w_{v(j)}}{l_n}\Big).
    }
The change in $Z_n(l)$ is not stationary, and decreases on the average
as $l$ increases, due to two reasons. First of all, the number of neutral
vertices decreases (as is apparent from the sum which is subtracted
in \refeq{Zl-asymp}), and the law of $\tilde w_{v(l)}$ becomes
stochastically smaller as $l$ increases. The latter can be understood by noting
that $\expec[\tilde w_{v(1)}]=(1+t n^{-1/3})\nu_n=1+tn^{-1/3}+o(n^{-1/3})$,
while $\frac{1}{n}\sum_{j\in[n]}\tilde w_{v(j)}=(1+t n^{-1/3})l_n/n$, and,
by Cauchy-Schwarz,
    \eqn{\label{csz1}
    l_n/n\approx \expec[W]\leq \expec[W^2]^{1/2} =\expec[W]^{1/2} \nu^{1/2}
    =\expec[W]^{1/2},
    }
so that $l_n/n\leq 1+o(1),$ and the inequality becomes strict when ${\rm Var}(W)>0$.
We now study these two effects in more detail.

The random variable $\Poi(\tilde w_{v(i)})-1$ has asymptotic mean
    \eqn{
    \expec[\Poi(\tilde w_{v(i)})-1]\sim \sum_{j\in [n]} \tilde w_j \prob(v(i)=j)-1
    =\sum_{j\in [n]} \tilde w_j\frac{w_j}{l_n}-1
    =\nu_n(1+tn^{-1/3})-1\approx 0.
    }
However, since we sum $\Theta(n^{2/3})$ contributions, and we multiply
by $n^{-1/3}$, we need to be rather precise and compute error terms up to
order $n^{-1/3}$ in the above computation. We shall do this rather precisely now,
by conditioning on $(v(j))_{j=1}^{i-1}$. Indeed,
    \eqan{
    \expec[\tilde w_{v(i)}-1]&\sim \nu_ntn^{-1/3}+\expec\big[\expec[w_{v(i)}-1\mid (v(j))_{j=1}^{i-1}]\big]\nn\\
    &\sim tn^{-1/3}+\expec\big[\sum_{l=1}^{n} w_l\indic{l\not\in \{v(j)\}_{j=1}^{i-1}}
    \frac{w_l}{l_n-\sum_{j=1}^{i-1} w_{v(j)}} \big]-1\nn\\
    &\sim tn^{-1/3}+\sum_{j\in [n]} \frac{w_j^2}{l_n}
    + \expec\big[\frac{1}{l_n^2}\sum_{j=1}^{i-1} w_{v(j)}\sum_{l=1}^{n} w_l\big] - \expec\big[\frac{1}{l_n}\sum_{j=1}^{i-1} w_{v(j)}^2\big]-1\nn\\
    &\sim tn^{-1/3}+i\big(\frac{\nu_n}{l_n}-\frac{1}{l_n}\sum_{j=1}^n w_j^3\big)
    \sim tn^{-1/3}+\frac{i}{l_n}\big(1-\frac{1}{l_n}\sum_{j=1}^n w_j^3\big).\label{exclusionn}
    }
When $i=\Theta(n^{2/3})$, these terms are indeed both
of order $n^{-1/3}$, and shall thus contribute to the scaling limit
of $(Z_n(l))_{l\geq 0}$.

The variance of $\Poi(\tilde w_{v(i)})$ is approximately equal to
    \eqan{
    {\rm Var}(\Poi(\tilde w_{v(i)}))
    &=\expec[{\rm Var}(\Poi(\tilde w_{v(i)}))\mid v(i)]+{\rm Var}(\expec[\Poi(\tilde w_{v(i)})\mid v(i)])\nn\\
    &=\expec[\tilde w_{v(i)}]+{\rm Var}(w_{v(i)})\sim \expec[\tilde w_{v(i)}^2]\sim \expec[w_{v(i)}^2],
    }
since $\expec[w_{v(i)}]=1+\Theta(n^{-1/3})$.
Summing the above over $i=1, \ldots, sn^{2/3}$ and multiplying by
$n^{-1/3}$ intuitively explains that
    \eqn{
    \lbeq{scaling-a}
    n^{-1/3}\sum_{i=1}^{sn^{2/3}} (\Poi(\tilde w_{v(i)})-1)\convd \sigma B(s) +st+\frac{s^2}{2\expec[W]}(1-\sigma^2),
    }
where we write $\sigma^2=\expec[W^3]/\expec[W]$ and we let $(B(s))_{s\geq 0}$
denote a standard Brownian motion. Note that when ${\rm Var}(W)>0$, then
$\expec[W]<1, \expec[W^2]>1,$ so that also $\expec[W^3]/\expec[W]>1$ and
the constant in front of $s^2$ is \emph{negative}.
We shall make the limit in \refeq{scaling-a} precise by
using a martingale functional central limit theorem.

The second term in
\refeq{Zl-asymp} turns out to be well-concentrated around its mean,
so that, in this heuristic, we shall replace it by its mean.
The concentration shall be proved using
concentration techniques on appropriate supermartingales.
This leads us to compute
    \eqan{
    \expec\Big[\sum_{i=1}^l \Poi\Big(\sum_{j=1}^{i-1} \frac{\tilde w_{v(i)}\tilde w_{v(j)}}{l_n}\Big)\Big]
    &\sim \expec\Big[\sum_{i=1}^l\sum_{j=1}^{i-1} \frac{\tilde w_{v(i)}\tilde w_{v(j)}}{l_n}\Big]
    \sim \expec\Big[\sum_{i=1}^l\sum_{j=1}^{i-1}\frac{w_{v(i)}w_{v(j)}}{l_n}\Big]\nn\\
    &\sim
    \frac 12 \expec\Big[\frac{1}{l_n}\Big(\sum_{j=1}^{i}w_{v(j)}\Big)^2\Big]
    \sim \frac{1}{2l_n} \expec\Big[\sum_{j=1}^{i}w_{v(j)}\Big]^2,
    }
the last asymptotic equality again following from the fact that the random
variable involved is concentrated.

We conclude that
    \eqn{
    \lbeq{scaling-b}
    n^{-1/3}\expec\Big[\sum_{i=1}^{sn^{2/3}} \Poi\Big(\sum_{j=1}^{i} \frac{\tilde w_{v(i)}\tilde w_{v(j)}}{l_n}\Big)\Big]
    \sim \frac{s^2}{2\expec[W]}.
    }
Subtracting \refeq{scaling-b} from \refeq{scaling-a}, these computations suggest, informally, that
    \eqn{
    \lbeq{infor-scaling}
    \bar{Z}_n(s) =  n^{-1/3} Z_n(\lfloor n^{2/3}s\rfloor)
    \convd \sigma B(s)+st-\frac{s^2\expec[W^3]}{2\expec[W]^2}
    =\sqrt{\frac{\expec[W^3]}{\expec[W]}}B(s)+st-\frac{s^2\expec[W^3]}{2\expec[W]^2},
    }
as required. Note the cancelation of the terms $\frac{s^2}{2\expec[W]}$
in \refeq{scaling-a} and \refeq{scaling-b}, where they appear with an opposite sign.
Our proof will make this analysis precise.

\subsection{Discussion}
\label{sec-disc}
Our results are generalizations of the critical behavior of Erd\H{o}s-R\'enyi random graphs,
which have received tremendous attention over the past decades. We refer
to \cite{Aldo97}, \cite{Boll01}, \cite{JanLucRuc00} and the references therein. Properties of the
limiting distribution of the largest component $\gamma_1(t)$ can be found in
\cite{Pitt01}, which, together with the recent local limit theorems in
\cite{HofKagMul09}, give excellent control over the joint tail behavior
of several of the largest connected components.

\paragraph{Comparison to results of Aldous.}
We have already discussed the relation between Theorem \ref{theo:main} and the results of Aldous on the largest connected components in the Erd\H{o}s-R\'enyi random graph. However, Theorem \ref{theo:main} is related to another result of Aldous \cite[Proposition 4]{Aldo97},
which is less well known, and which investigates a kind of Norros-Reittu model (see \cite{NorRei06})
for which the ordered \emph{weights} of the clusters are determined. Here,
the weight of a set of vertices $A\subseteq [n]$ is defined by
$\bar{w}_{\sss A}=\sum_{a\in A} w_a$. Indeed, Aldous defines an inhomogeneous
random graph where the edge probability is equal to
    \eqn{
    p_{ij}=1-{\mathrm e}^{-q x_ix_j},
    }
and assumes that the pair $(q,(x_i)_{i=1}^n)$ satisfies the following scaling relation:
    \eqn{
    \lbeq{Aldo-defs}
    \frac{\sum_{i=1}^n x_i^3}{\big(\sum_{i=1}^n x_i^2\big)^3}
    \rightarrow 1,
    \qquad
    q-\Big(\sum_{i=1}^n x_i^2\Big)^{-1}\rightarrow t,
    \qquad
    \max_{j\in [n]}x_j=o\Big(\sum_{i=1}^n x_i^2\Big).
    }
When we pick
    \eqn{
    x_j=w_j \frac{\big(\sum_{i=1}^n w_i^3\big)^{1/3}}{\sum_{i=1}^n w_i^2},
    \qquad
    q=\frac{\big(\sum_{i=1}^n w_i^2\big)^{2}}{\big(\sum_{i=1}^n w_i^3\big)^{2/3}l_n}(1+tn^{-1/3}),
    }
then these assumptions are very similar to conditions (a)-(c).
However, the asymptotics of $q$ in \refeq{Aldo-defs}
is replaced with
    \eqn{
    q-\Big(\sum_{i=1}^n x_i^2\Big)^{-1}
    =\frac{\frac{1}{n}\sum_{i=1}^n w_i^2}{\big(\frac{1}{n}\sum_{i=1}^n w_i^3\big)^{2/3}}(n^{1/3}\nu_n(1+tn^{-1/3})-n^{1/3})
    \rightarrow \frac{\expec[W^2]}{\expec[W^3]^{2/3}}t=\frac{\expec[W]}{\expec[W^3]^{2/3}}t,
    }
where the last equality follows from the fact that $\nu=\expec[W^2]/\expec[W]=1$.
This scaling in $t$ simply means that the parameter $t$ in the
process $W^t_*(s)$ in \eqref{eqn:wt} is rescaled, which is explained in
more detail in the scaling relations in \refeq{eq-scal}. Write $\CC_n^i(t)$
for the component with the $i^{\rm th}$ largest {\bf weight}, and let
$\bar{w}_{\sss\CC_n^i(t)}=\sum_{j\in \CC_n^i(t)}w_j$
denote the cluster weight. Then, Aldous \cite[Proposition 4]{Aldo97} proves that
    \eqn{
    \Big(\frac{\big(\sum_{i=1}^n w_i^3\big)^{1/3}}{\sum_{i=1}^n w_i^2}\bar{w}_{\sss\CC_n^i(t)}\Big)_{i\geq 1}
    \convd \big(\gamma_i(t\expec[W]/\expec[W^3]^{2/3})\big)_{i\geq 1},
    }
where we recall
that $\big(\gamma_i(t)\big)_{\i\geq 1}$ is the scaling limit of the ordered component sizes
in the Erd\H{o}s-R\'enyi random graph with parameter $p=(1+tn^{-1/3})/n$. Now,
    \eqn{
    \lbeq{scaling-Aldous}
    \frac{\big(\sum_{i=1}^n w_i^3\big)^{1/3}}{\sum_{i=1}^n w_i^2}
    \sim n^{-2/3} \expec[W^3]^{1/3}/\expec[W^2]=n^{2/3} \expec[W^3]^{1/3}/\expec[W],
    }
and one would expect that $\bar{w}_{\sss\CC_n^i(t)}\sim \CC_n^i(t)$, which is consistent
with \refeq{infor-scaling} and \refeq{eq-scal}.

\paragraph{Related models.}
The model studied here is asymptotically equivalent to many related models
appearing in the literature,
for example to the  \emph{random graph with
prescribed expected degrees} that has been studied intensively by Chung and Lu (see
\cite{ChuLu02a,ChuLu02b,ChuLu03, ChuLu06c, ChuLu06}).
This model corresponds to the rank-1 case of the general inhomogeneous random graphs
studied in \cite{BolJanRio07}. Here
    \eqn{
        \lbeq{pij-NR}
 p_{ij} =\min\left\{\frac{w_i w_j}{l_n},1\right\},
    }
and the \emph{generalized random graph}
\cite{BriDeiMar-Lof05}, for which
    \eqn{
    \lbeq{pij-NR}
    p_{ij} = \frac{w_i w_j}{l_n+w_iw_j},
    }
See \cite[Section 2]{Hofs09a}, which in turn is based on
\cite{Jans08a}, for more details
on the asymptotic equivalence of such inhomogeneous random graphs.
Further, Nachmias and Peres \cite{NacPer07a} recently proved similar scaling limits for
critical percolation on random regular graphs.

\paragraph{Alternative approach by Turova.}
Turova \cite{Turo09} recently obtained results for a setting that is similar to ours. Turova takes
the edge probabilities to be $p_{ij}=\min\{x_ix_j/n,1\}$, and assumes that
$(x_i)_{i=1}^n$ are i.i.d.\ random variables with $\expec[X^3]<\infty$. This setting follows from ours
by taking
    \eqn{
    w_i=x_i\big(\frac 1n\sum_{j=1}^n x_j\big)^{1/2}.
    }
Naturally, the critical point changes in Turova's setting, and is equal to $\expec[X^2]=1$.

First versions of the paper \cite{Turo09} and this paper were uploaded almost simultaneously on the ArXiv.
Comparing the two papers gives interesting insights in how to deal with the inherent size-biased orderings in two rather different ways.
Turova applies discrete martingale techniques in the spirit of Martin-L\"{o}f's \cite{Mart98} work on diffusion approximations for critical epidemics, while our approach is more along the lines of the original paper of Aldous \cite{Aldo97}, relying on concentration techniques and supermartingales (see Lemma \eqref{lem-wv(i)-sums}).  Further, our result is
slightly more general that the one in \cite{Turo09}. In fact, our discussions with Turova inspired us to extend our setting to one that includes i.i.d.~weights (Turova's original setting). We should also mention that Turova's first identification of the scaling limit was missing a factor $\expec[X^3]$ in the drift term (which was corrected in a later version). This factor arises from rather subtle effects of the sized-biased orderings.

\paragraph{The necessity of conditions (a)-(c).} The conditions (a)-(c) 
provide conditions under which we prove convergence. One may wonder whether
these conditions are merely sufficient, or also necessary. Condition (b) 
gives stability of the weight structure, which implies that
the local neighborhoods in our random graphs locally converge to appropriate branching
processes. The latter is a strengthening of the assumption that our random
graphs are sparse, and is a natural condition to start with. We believe that, given that
condition (b) holds, conditions (c) and (a) are necessary. Indeed, 
Aldous and Limic give several examples where the scaling of the largest critical
cluster is $n^{2/3}$ with a \emph{different} scaling limit when $w_1n^{1/3}\rightarrow c_1>0$
(see \cite[Proof of lemma 8, p. 10]{aldous-limic}).
Therefore, for Theorem \ref{theo:main} to hold (with the prescribed scaling limit in terms
of ordered Brownian excursions), condition (a) seems to be necessary. Since conditions
(b) and (c) imply condition (a), it follows that if we assume condition (b), then
we need the other two conditions for our main result to hold. This answers 
\cite[Open problem (2), p.\ 851]{Aldo97}.

\paragraph{Inhomogeneous random graphs with infinite third moments.}
In the present setting, when it is assumed that $\expec[W^3]<\infty$, the scaling
limit turns out to be a scaled version of the scaling limit for the
Erd\H{o}s-R\'enyi random graph as identified in \cite{Aldo97}.
In \cite{BhaHofLee09b}, we have recently studied the case where
$\expec[W^3]=\infty$, for which the critical behavior
turns out to be fundamentally different. Indeed, when $W$ has a power 
law with exponent $\tau\in(3,4)$, the clusters have asymptotic 
size $n^{\frac{\tau-2}{\tau-1}}$ (see \cite{Hofs09a}).
The scaling limit itself turns out to be a so-called `thinned' L\'{e}vy
process, that consists of infinitely many Poisson processes of which
only the first event is counted, which already appeared in 
\cite{aldous-limic} in the context of random graphs having $n^{2/3}$
critical behavior. Moreover, we prove in \cite{BhaHofLee09b}
that the vertex $i$ is in the largest connected component with non-vanishing
probability as $n\rightarrow \infty$, which implies that
the highest weight vertices characterize the largest components
(`power to the wealthy'). This is in
sharp contrast to the present setting, where the probability that vertex 1
(with the largest weight) is in the largest component is negligible, and
instead the largest connected component is an extreme value event arising
from many trials with roughly equal probability (`power to the masses').

\section{Weak convergence of cluster exploration}
\label{sec-cluster-explor-tau>4}
In this section, we shall study the scaling limit of the cluster
exploration studied in Section \ref{sec-cluster-exploration} above.
The main result in this paper is the following theorem:
\begin{Theorem}[Weak convergence of cluster exploration]
\label{theo:walk-bm}
Assume that the weight sequence $\bfw$
satisfies conditions {\rm(a)}, {\rm(b)} and {\rm(c)}. Consider the breadth-first walk $Z_n(\cdot)$ of \eqref{eqn:walk-comp}
exploring the components of the random graph $\GG_n^t(\bfw)$. Define
	\eqn{
    \bar{Z}_n(s) = n^{-1/3} Z_n(\lfloor n^{2/3}s\rfloor).
    }
Then, as $n\to\infty$,
	\eqn{
    \bar{Z}_n \convd W^t_*,
    }
where $W^t_*$ is the process defined in \eqref{eqn:wt},
in the sense of convergence in the $J_1$ Skorohod topology on the
space of right-continuous left-limited functions on ${\mathbb R}^+$.
\end{Theorem}

Assume this theorem for the time being and let us show how
this immediately proves Theorem \ref{theo:main}. Comparing
\eqref{eqn:bt} and \eqref{eqn:walk-comp}, Theorem \ref{theo:walk-bm}
suggests that also the excursions of $\bar{Z}_n$ beyond
past minima arranged in increasing order converge to the corresponding
excursions of $W^t_*$ beyond past minima arranged in increasing order.
See Aldous \cite[Section 3.3]{Aldo97} for a proof of this fact.
Therefore, Theorem \ref{theo:walk-bm} implies Theorem \ref{theo:main}.
The remainder of this paper is devoted to the proof of Theorem \ref{theo:walk-bm}.

\vspace{.2cm}
\noindent
{\it Proof of Theorem \ref{theo:walk-bm}.} We shall make use of a martingale
central limit theorem. From Equation \eqref{eqn:pij-pijst} we had
    \eqn{
    \lbeq{pproxx}
    p_{ij} \approx \left(1+\frac{t}{n^{1/3}}\right)\frac{w_i w_j}{l_n},
    }
and we shall use the above as an equality for the rest of the proof as this shall simplify exposition. It is quite easy to show that the error made is negligible in the limit.

 Recall from \refeq{Zn-def} that
    \eqn{
    Z_n(k) = \sum_{i=0}^k (c(i)-1).
    }
Then, we decompose
    \eqn{
    \lbeq{Zk-split}
    Z_n(k)=M_n(k)+A_n(k),
    }
where
    \eqn{
    \lbeq{Ak-def}
    M_n(k)=\sum_{i=0}^k (c(i)-\expec[c(i)\mid {\cal F}_{i-1}]),
    \qquad
    A_n(k)=\sum_{i=0}^k \expec[c(i)-1\mid {\cal F}_{i-1}],
    }
with ${\cal F}_{i}$ the natural filtration of $Z_n$.
Then, clearly, $\{M_n(k)\}_{k=0}^n$ is a martingale. For a process
$\{S_k\}_{k=0}^n$, we further write
    \eqn{
    \bar{S}_n(u)=n^{-1/3} S_n(\lfloor un^{2/3} \rfloor).
    }
Furthermore, let
    \eqn{
    \lbeq{Bk-def}
    B_n(k)=\sum_{i=0}^k \left(\expec[c(i)^2\mid {\cal F}_{i-1}]-\expec[c(i)\mid {\cal F}_{i-1}]^2\right).
    }
Then, by the martingale central limit theorem (\cite[Theorem 7.1.4]{EthKur86}),
Theorem \ref{theo:walk-bm} follows when the following three conditions hold:
    \eqan{
    \sup_{s\leq u} \Big|\bar{A}_n(s) + \frac{s^2\sigma_3}{2\mu^2} - st\Big| &\convp 0,
    \label{eqn:an-to-0}\\
	n^{-2/3}B_n(n^{2/3}u) &\convp \frac{\sigma_3u}{\mu}, \label{eqn:bn-u}\\
	\expec(\sup_{s\leq u} |\bar{M}_n(s) - \bar{M}_n(s-)|^2) &\to 0. \label{eqn:mn-mn}
    }
Indeed, the last two equations, by \cite[Theorem 7.1.4]{EthKur86} imply
that the process $\bar{M}_n(s) = n^{-1/3}M_n(n^{2/3 }s)$ satisfies the asymptotics
    \eqn{
    \bar{M}_n\convd \sqrt{\frac{\sigma_3}{\mu}}B,
    }
where as before $B$ is standard Brownian motion, while \eqref{eqn:an-to-0} gives the drift term in
\eqref{eqn:wt} and this completes the proof.

We shall now start to verify the conditions \eqref{eqn:an-to-0}, \eqref{eqn:bn-u} and
\eqref{eqn:mn-mn}. Throughout the proof, we shall assume, without loss of generality, that
$w_1\geq w_2\geq  \ldots\geq w_n$. Recall that we shall work with
weight sequence $\tilde \bfw=(1+t n^{-1/3})\bfw$, for which the edge
probabilities are approximately equal to $w_iw_j(1+t n^{-1/3})/l_n$ (recall \eqref{pproxx}).

We note that, since $M_n(k)$ is a discrete martingale,
    \eqan{
    \sup_{s\leq u} |\bar{M}_n(s) - \bar{M}_n(s-)|^2&=n^{-2/3} \sup_{k\leq un^{2/3}}
    (M_n(k)-M_n(k-1))^2
    \leq n^{-2/3} (1+\sup_{k\leq un^{2/3}}c(k)^2)\nn\\
    &\leq n^{-2/3} (1+\Delta_n^2),
    }
where $\Delta_n$ is the maximal degree in the graph. It is not hard to see that,
by condition (a), $\tilde w_i=o(n^{1/3})$, so that
    \eqn{
    \expec(\sup_{s\leq u} |\bar{M}_n(s) - \bar{M}_n(s-)|^2)
    \leq n^{-2/3}(1+\expec[\Delta_n^2])=o(n^{-2/3} n^{2/3})=o(1).
    }
This proves \eqref{eqn:mn-mn}.

We continue with \eqref{eqn:an-to-0} and \eqref{eqn:bn-u}, for which we first analyse
$c(i)$. In the course of the proof, we shall make use of the following lemma, which
lies at the core of the argument:

\begin{Lemma}[Sums over sized-biased orderings]
\label{lem-wv(i)-sums}
As $n\rightarrow \infty$, for all $u>0$,
    \eqan{
    \sup_{u\leq t}\big|n^{-2/3}\sum_{i=1}^{n^{2/3}u}w_{v(i)}^2-\frac{\sigma_3 u}{\mu}|&\convp 0,\label{sum-1-conv}\\
    n^{-2/3}\sum_{i=1}^{n^{2/3}u} \expec[w_{v(i)}^2\mid {\cal F}_{i-1}]&\convp
    \frac{\sigma_3 u}{\mu}. \label{sum-2-conv}
    }
\end{Lemma}

\proof We start by proving \eqref{sum-1-conv}, for which we write
    \eqn{
    H_n(u)= n^{-2/3}\sum_{i=1}^{\lfloor u n^{2/3} \rfloor} w_{v(i)}^2.
    }
We shall use a randomization trick introduced by Aldous \cite{Aldo97}.
Indeed, let $T_j$ be a sequence of independent exponential
random variables with rate $w_j/l_n$ and define
	\eqn{
    \tilde{H}_n(v) =  n^{-2/3}\sum_{i=1}^n w_j^2 \ind\{T_j \leq n^{2/3}v\}.
    }
Note that by the properties of the exponential random variables,
if we rank the vertices according to the order in which they arrive then
they appear in size-biased order. More precisely, for any $v$,
	\eqn{
    \sum_{j=1}^n w_j^2 \ind\{T_j \leq n^{2/3}v\} = \sum_{i=1}^{N(v n^{2/3})} w_{v(i)}^2=
    H_n(N(v n^{2/3})),
    }
where
	\eqn{
    \lbeq{Nt-def}
    N(t):= \#\{j: T_j \leq t\}.
    }
As a result, when $N(2tn^{2/3})\geq tn^{2/3}$ {\bf whp}, we have that
    \eqan{
    \sup_{u\leq t}\big|n^{-2/3}\sum_{i=1}^{n^{2/3}u}w_{v(i)}^2-\sigma_3 u|
    &\leq \sup_{u\leq 2t}\big|n^{-2/3}\sum_{i=1}^{N(u n^{2/3})}w_{v(i)}^2-\sigma_3 n^{-2/3}N(u n^{2/3})|\nn\\
    &\leq \sup_{u\leq 2t}\big|n^{-2/3}\tilde{H}_n(u)-\sigma_3 u|
    +\sigma_3\sup_{u\leq 2t}\big|n^{-2/3}N(u n^{2/3})-u|.
    }
We shall prove that both terms converge to zero in probability. We start with the second,
for which we use that the process
	\eqn{
    Y_0(s) = \frac{1}{n^{1/3}} \left(N(s n^{2/3}) - s n^{2/3}\right),
    }
is a supermartingale, since
    \eqan{
    \expec[N(t+s)\mid {\cal F}_t]&=
    N(t)+\expec[N(t+s)-N(t)\mid {\cal F}_t]
    \leq \expec[\#\{j: T_j \in(t, t+s]\}\mid {\cal F}_t]\nn\\
    &\leq\sum_{j=1}^{n} (1-{\mathrm e}^{-w_js/l_n})
    \leq \sum_{j=1}^{n} \frac{w_j s}{l_n}=s,
    }
as required. Therefore,
	\eqn{
    |\expec[Y_0(t)]| = - \expec[Y_0(t)] =\frac{1}{n^{1/3}}\left[ tn^{2/3}-\sum_{i=1}^{n} (1-\exp(-t n^{2/3} w_i/l_n))\right].
    }
Using the fact that $1-{\mathrm e}^{-x} \leq x - x^2/2,$
we obtain that, also using that $\nu_n=1+o(1)$,
	\eqn{
    |\expec[Y_0(t)]| \leq \sum_{i=1}^n \frac{w_i^2t^2}{2l_n^2}
    =\frac{n\nu_n}{l_n} \frac{t^2}{2}=\frac{t^2}{2\mu}+o(1).
    }
Similarly, by the independence of $\{T_j\}_{j\in [n]}$,
	\eqan{
    \mbox{Var}(Y_0(t))&= n^{-2/3}\mbox{Var}(N(s n^{2/3}))=
    n^{-2/3}\sum_{j=1}^n \prob(T_j\leq tn^{2/3})(1-\prob(T_j\leq tn^{2/3}))\nn\\
    &\leq n^{-2/3}\sum_{j=1}^n \frac{w_jtn^{2/3}}{l_n}=t.
    }
Now we use the supermartingale inequality (Aldous \cite[page 831, proof of Lemma 12]{Aldo97}),
stating that, for any supermartingale $Y=(Y(s))_{s\geq 0}$,
	\eqn{\eps \prob(\sup_{s\leq t} |Y(s)|> 3\eps) \leq 3 \expec(|Y(t)|) \leq 3 \left(|\expec(Y(t))| +
    \sqrt{\mbox{Var}(Y(t))}\right).
    \label{eqn:supermart-ineq}
    }
Equation \eqref{eqn:supermart-ineq} shows that, for any large $A$,
	\eqn{
    \prob(\sup_{s\leq t}|N(s n^{2/3}) - sn^{2/3}| > 3A n^{1/3}) \leq \frac{3 (t^2/{2\mu} + t)}{A}.
    }
This clearly proves that, for every $t>0$,
    \eqn{
    \lbeq{convp-Ns}
    \sup_{u\leq 2t}\big|n^{-2/3}N(u n^{2/3})-u|\convp 0.
    }
Observe that \refeq{convp-Ns} also immediately proves that,
{\bf whp}, $N(2tn^{2/3})\geq tn^{2/3}$.

To deal with $\tilde{H}_n(v)$, we define
	\eqn{
    Y_1(u) = \tilde{H}_n(u) - \mu_3(n) u,
    }
where
    \eqn{
    \mu_3(n)=\sum_{j=1}^n \frac{w_j^3}{l_n}=\frac{\sigma_3}{\mu}+o(1),
    }
and note that $Y_1(u)$ is a supermartingale. Indeed,
writing $\FF_t$ to be the natural filtration of the above process, we have,
for $s<t$ and letting $V_s = \set{v: T_v< s n^{2/3}}$
    \eqn{
    \expec(Y_1(t)|\FF_s)= Y_1(s) + \frac{1}{n^{2/3}}
    \sum_{j\notin V_s} w_j^2 \left(1-\exp\left(-\frac{n^{2/3}(t-s) w_j}{l_n}\right)\right) - \mu_3(n)(t-s).
    }
Now using the inequality $1-{\mathrm e}^{-x} \leq x$ for $x\in[0,1]$ we get that
    \eqn{
    \expec(Y_1(t)|\FF_s)\leq Y_1(s),
    }
as required. Again we can easily compute, using condition (a), that
    \eqan{
    |\expec[Y_1(t)]|&=- \expec[Y_1(t)] =
    \mu_3(n) t-n^{-2/3}\sum_{i=1}^{n} w_i^2(1-\exp(-t n^{2/3} w_i/l_n))\nn\\
    &=n^{-2/3}\sum_{i=1}^{n} w_i^2\big(\exp(-t n^{2/3} w_i/l_n)-1+tw_i\big)\nn\\
    &\leq n^{-2/3}\sum_{i=1}^{n} w_i^2 \frac{(t n^{2/3} w_i)^2}{2l_n^2}
    \leq n^{2/3} \sum_{i=1}^n \frac{w_i^4}{2l_n^2}=o(n^{2/3} n^{1/3})\sum_{i=1}^n \frac{w_i^3}{l_n^2}
    =o(1).
    }
By independence,
    \eqan{
    \mbox{Var}(Y_1(t))&=n^{-4/3}\sum_{j=1}^n w_i^4 (1-\exp(-t n^{2/3} w_i/l_n))\exp(-t n^{2/3} w_i/l_n)\nn\\
    &\leq n^{-2/3} t\sum_{j=1}^n \frac{w_i^5}{l_n}=o(1) \sum_{j=1}^n \frac{w_i^3}{l_n}
    =o(1).
    }
Therefore, \eqref{eqn:supermart-ineq} completes the proof of \eqref{sum-1-conv}.

The proof of \eqref{sum-2-conv} is a little easier. We denote
    \eqn{
    {\cal V}_{i}=\{v(j)\}_{j=1}^i.
    }
Then, we compute explicitly
    \eqan{
    \expec[w_{v(i)}^2\mid {\cal F}_{i-1}]&=\sum_{j\in[n]} w_{j}^2\prob(v(i)=j\mid {\cal F}_{i-1})\\
    &=\frac{\sum_{j\not\in {\cal V}_{i-1}} w_{j}^3}{\sum_{j\not\in {\cal V}_{i-1}}w_j}.
    }
Now, uniformly in $i\leq sn^{2/3}$, again using condition (a),
    \eqn{
    \sum_{j\not\in {\cal V}_{i-1}}w_j=\sum_{j\in[n]}w_j+O((\max_{j\in [n]} w_j) i)
    =l_n+o(n)
    }
for every $i\leq sn^{2/3}$ and since $1/(\tau-1)<1/3$ for $\tau>4$. Similarly,
again uniformly in $i\leq sn^{2/3}$, and using that $j\mapsto w_j$ is non-increasing,
    \eqn{
    \lbeq{bound-w3}
    \Big|\sum_{j\not\in {\cal V}_{i-1}}w_j^3-l_n\sigma_3(n)\Big|
    \leq \sum_{j=1}^{sn^{2/3}} w_j^3
    =o(n).
    }
Indeed, we shall show that conditions (b) and (c) imply that
    \eqn{
    \lbeq{unif-int-w3}
    \lim_{K\rightarrow \infty} \lim_{n\rightarrow \infty} \frac{1}{n}\sum_{j\in [n]} \indic{w_j>K}w_j^3=0.
    }
Equation \refeq{unif-int-w3} in particular implies that $w_1^3=o(n)$, so that conditions (b) and (c) also
imply condition (a). By the weak convergence stated in condition (b),
    \eqn{
    \frac{1}{n}\sum_{j\in [n]} \indic{w_j\leq K} w_j^3=\expec[\indic{w_{V_n}\leq K} w_{V_n}^3]
    \rightarrow \expec[\indic{W\leq K} W^3].
    }
As a result, we have that
    \eqn{
    \lim_{n\rightarrow \infty} \frac{1}{n}\sum_{j\in [n]} \indic{w_j>K}w_j^3=
    \lim_{n\rightarrow \infty} \frac{1}{n}\sum_{j\in [n]}w_j^3-
    \lim_{n\rightarrow \infty}\frac{1}{n}\sum_{j\in [n]} \indic{w_j\leq K}w_j^3
    =\expec[\indic{W>K} W^3].
    }
The latter converges to 0 when $K\rightarrow \infty$, since $\expec[W^3]<\infty$.
We finally show that \refeq{unif-int-w3} implies that $\sum_{j=1}^{sn^{2/3}} w_j^3
=o(n)$. For this, we note that, for each $K$,
    \eqn{
    \frac{1}{n}\sum_{j=1}^{sn^{2/3}} w_j^3\leq
      \frac{1}{n}\sum_{j=1}^{sn^{2/3}} \indic{w_j\leq K}w_j^3
   + \frac{1}{n}\sum_{j\in [n]} \indic{w_j>K} w_j^3\leq    K^3 sn^{-1/3} +\frac{1}{n}\sum_{j\in [n]} \indic{w_j>K} w_j^3
 =o(1),
    }
when we first let $n\rightarrow\infty$, followed by $K\rightarrow \infty$.

We conclude that, uniformly for $i\leq sn^{2/3}$,
    \eqn{
    \expec[w_{v(i)}^2\mid {\cal F}_{i-1}]=\frac{\sigma_3}{\mu}+\op(1).
    }
This proves \eqref{sum-2-conv}.
\qed

To complete the proof of Theorem \ref{theo:walk-bm},
we proceed to investigate $c(i)$. By construction,
we have that, conditionally on ${\cal V}_{i}$,
    \eqn{
    c(i)\stackrel{d}{=}\sum_{j\not\in {\cal V}_{i}} I_{ij},
    }
where $I_{ij}$ are (conditionally) independent indicators with
    \eqn{
    \prob(I_{ij}=1\mid {\cal V}_{i})=\frac{w_{v(i)}w_j(1+t n^{-1/3})}{l_n},
    }
for all $j\not\in {\cal V}_i$. Furthermore, when we condition on
${\cal F}_{i-1}$, we know ${\cal V}_{i-1}$, and we have that,
for all $j\not\in {\cal V}_{i-1}$,
    \eqn{
    \prob(v(i)=j\mid {\cal F}_{i-1})=\frac{w_j}{\sum_{s\in{\cal V}_{i-1}^c} w_s}.
    }
Since ${\cal V}_i={\cal V}_{i-1}\cup \{v(i)\}={\cal V}_{i-1}\cup \{j\}$
when $v(i)=j$, this gives us all we need to
know to compute conditional expectations involving
$c(i)$ given ${\cal F}_{i-1}$.

Now we start to prove \eqref{eqn:an-to-0}, for which we note that
    \eqan{
    \expec[c(i)\mid {\cal F}_{i-1}]
    &=\sum_{j\in {\cal V}_{i-1}^c} \prob(v(i)=j\mid {\cal F}_{i-1}) \expec[c(i)\mid {\cal F}_{i-1}, v(i)=j]\nn\\
    &=\sum_{j\in {\cal V}_{i-1}^c} \prob(v(i)=j\mid {\cal F}_{i-1}) \sum_{l\not\in{\cal V}_{i-1}\cup \{j\}}\frac{\tilde w_{j}w_l}{l_n}.
    }
Then we split
    \eqan{
    \lbeq{c(i)-1-cond-expec}
    \expec[c(i)-1\mid {\cal F}_{i-1}]
    &=\sum_{j\in {\cal V}_{i-1}^c} \prob(v(i)=j\mid {\cal F}_{i-1})\tilde w_j-1
    -\sum_{j\in {\cal V}_{i-1}^c} \prob(v(i)=j\mid {\cal F}_{i-1})\tilde w_j\sum_{l\in {\cal V}_{i-1}\cup \{j\}}\frac{w_l}{l_n}\\
    &=\expec[\tilde w_{v(i)}-1\mid {\cal F}_{i-1}]-\expec[\tilde w_{v(i)}\mid {\cal F}_{i-1}]\sum_{s=1}^{i-1} \frac{w_{v(s)}}{l_n}
    -\expec[\frac{w_{v(i)}^2(1+t n^{-1/3})}{l_n}\mid {\cal F}_{i-1}].\nn
    }
By condition (a), the last term is bounded by
$\Op(w_1^2/l_n)=\op(n^{-1/3})$  and is therefore an error term.
We continue to compute
    \eqan{
    \lbeq{wv(i)-comp}
    \expec[\tilde w_{v(i)}-1\mid {\cal F}_{i-1}]
    &=\sum_{j\in {\cal V}_{i-1}^c} \frac{w_j^2(1+t n^{-1/3})}{\sum_{s\in {\cal V}_{i-1}^c} w_s}-1\\
    &=\sum_{j\in {\cal V}_{i-1}^c} \frac{w_j^2(1+t n^{-1/3})}{l_n}-1
    +\sum_{j\in {\cal V}_{i-1}^c} \frac{w_j^2(1+t n^{-1/3})}{l_n\sum_{s\in {\cal V}_{i-1}^c} w_s}\Big(\sum_{s\in {\cal V}_{i-1}} w_s\Big).\nn
    }
The last term equals
    \eqn{
    \expec[\frac{\tilde w_{v(i)}}{l_n}\mid {\cal F}_{i-1}]\sum_{s=1}^{i-1} w_{v(s)},
    }
which equals the second term in \refeq{c(i)-1-cond-expec}, and thus these two contributions
cancel in \refeq{c(i)-1-cond-expec}. This exact cancelation is in the spirit of
the one discussed below \refeq{infor-scaling}. Therefore,
writing $\tilde\nu_n=\nu_n(1+t n^{-1/3})$,
    \eqan{
    \lbeq{c(i)-1-cond-expec-cont}
    \expec[c(i)-1\mid {\cal F}_{i-1}]
    &=\sum_{j\in {\cal V}_{i-1}^c} \frac{w_j^2(1+t n^{-1/3})}{l_n}-1+\op(n^{-1/3})\nn\\
    &=\sum_{j=1}^n \frac{w_j^2(1+t n^{-1/3})}{l_n}-1 -\sum_{j\in {\cal V}_{i-1}}
    \frac{w_j^2(1+t n^{-1/3})}{l_n} +\op(n^{-1/3})\nn\\
    &=(\tilde\nu_n-1) -\sum_{s=1}^{i-1} \frac{w_{v(s)}^2(1+t n^{-1/3})}{l_n}+\op(n^{-1/3})\nn\\
    &=(\tilde\nu_n-1) -\sum_{s=1}^{i-1} \frac{w_{v(s)}^2}{l_n}+\op(n^{-1/3}).
    }
As a result, we obtain that
    \eqan{
    A_n(k)&=\sum_{i=0}^k \expec[c(i)-1\mid {\cal F}_{i-1}]
    =k(\tilde\nu_n-1)-\sum_{i=0}^k \sum_{s=1}^{i-1} \frac{w_{v(s)}^2}{l_n} +\op(k n^{-1/3}).
    }
Thus,
    \eqan{
    \bar{A}_n(s)&=ts -n^{-1/3}\sum_{i=0}^{sn^{2/3}} \sum_{l=1}^{i-1} \frac{w_{v(l)}^2}{l_n}+\op(1).
    }
By \eqref{sum-1-conv} in Lemma \ref{lem-wv(i)-sums}, we have that
    \eqn{
    \sup_{t\leq u}|n^{-2/3}\sum_{s=1}^{tn^{2/3}} w_{v(s)}^2-\frac{\sigma_3}{\mu} t|\convp 0,
    }
so that
    \eqn{
    \sup_{t\leq u}\big|\bar{A}_n(s)-ts +\frac{s^2 \sigma_3}{2\mu^2}\big|\convp 0.
    }
This proves \eqref{eqn:an-to-0}.

The proof for \eqref{eqn:bn-u} is similar, and we start by
noting that \refeq{c(i)-1-cond-expec-cont} gives that
    \eqn{
    B_n(k)=\sum_{i=0}^k \expec[c(i)^2\mid {\cal F}_{i-1}]-\expec[c(i)\mid {\cal F}_{i-1}]^2
    =\sum_{i=0}^k \expec[c(i)^2-1\mid {\cal F}_{i-1}]+\Op(kn^{-1/3}).
    }
Now, as above, we obtain that
    \eqan{
    \expec[c(i)^2\mid {\cal F}_{i-1}]
    &=\sum_{j\in {\cal V}_{i-1}^c} \prob(v(i)=j\mid {\cal F}_{i-1})
    \sum_{\stackrel{s_1, s_2\not\in{\cal V}_{i-1}\cup \{j\}}{s_1\neq s_2}}\frac{\tilde w_{j}w_{s_1}}{l_n}\frac{\tilde w_{j}w_{s_2}}{l_n}\\
    &\qquad+\sum_{j\in {\cal V}_{i-1}^c} \prob(v(i)=j\mid {\cal F}_{i-1})
    \sum_{s\not\in{\cal V}_{i-1}\cup \{j\}}\frac{\tilde w_{j}w_{s}}{l_n}.\nn
    }
We compute the second term as, using condition (a),
    \eqn{
    \sum_{j\in {\cal V}_{i-1}^c} \prob(v(i)=j\mid {\cal F}_{i-1})
    \sum_{s\not\in{\cal V}_{i-1}\cup \{j\}}\frac{\tilde w_{j}w_{s}}{l_n}
    =\sum_{j\in {\cal V}_{i-1}^c} \prob(v(i)=j\mid {\cal F}_{i-1})w_j(1+n\op(in^{-2/3}))
    =1+\op(1).
    }
The first sum is similarly computed as
    \eqan{
    \sum_{j\in {\cal V}_{i-1}^c} \prob(v(i)&=j\mid {\cal F}_{i-1})
    \sum_{\stackrel{s_1, s_2\not\in{\cal V}_{i-1}\cup \{j\}}{s_1\neq s_2}}\frac{\tilde w_{j}w_{s_1}}{l_n}\frac{\tilde w_{j}w_{s_2}}{l_n}\\
    &=\sum_{j\in {\cal V}_{i-1}^c} \prob(v(i)=j\mid {\cal F}_{i-1}) \tilde w_j^2 (1+\op(1))
    =\expec[\tilde w_{v(i)}^2\mid {\cal F}_{i-1}]+\op(1),\nn
    }
so that
    \eqan{
    n^{-2/3}B_n(n^{2/3}u)
    &=n^{-2/3}\sum_{i=1}^{n^{2/3}u} \expec[w_{v(i)}^2\mid {\cal F}_{i-1}]+\op(1)=\frac{\sigma_3}{\mu} u+\op(1),
    }
where the last equality follows from \eqref{sum-2-conv} in Lemma \ref{lem-wv(i)-sums}.
The proofs of \eqref{eqn:an-to-0}, \eqref{eqn:bn-u} and \eqref{eqn:mn-mn} complete the proof of
Theorem \ref{theo:walk-bm}.
\qed

\section{Verification of conditions (b)--(c): Proof of Corollary \ref{cor:main}}
\label{sec-ver-cond}
\subsection{Verification of conditions for i.i.d.\ weights}
We now check conditions (b) and (c) for the case that
$\bfw=(W_1, \ldots, W_n)$ where $\{W_i\}_{i\in [n]}$ are i.i.d.\ random variables
with $\expec[W^3]<\infty$.
Condition (b) is equivalent to the a.s.\ convergence of the empirical distribution function,
while \refeq{conv-third-mom} in condition (c) holds by the strong law of large numbers.
Equation \refeq{conv-first-mom} in condition (c) holds by the central limit theorem
(even with $\op(n^{-1/3})$ replaced with $\Op(n^{-1/2}))$. The bound in
\refeq{conv-sec-mom} in condition (c) is a bit more delicate.

To bound \refeq{conv-sec-mom} in condition (c), we will first split
    \eqan{
    \frac{1}{n}\sum_{i=1}^n \big(W_i^2-\expec[W^2])&=\frac{1}{n}\sum_{i=1}^n \Big(W_i^2\indic{W_i\leq n^{1/3}}-\expec[W^2\indic{W\leq n^{1/3}}]\Big)\\
    &\qquad+\frac{1}{n}\sum_{i=1}^n \Big(W_i^2\indic{W_i>n^{1/3}}-\expec[W^2\indic{W>n^{1/3}}]\Big).\nn
    }
On the second term, we use the Markov inequality to show that,
for every $\vep>0$,
    \eqan{
    &\prob\Big(\frac{1}{n}\sum_{i=1}^n \Big|W_i^2\indic{W_i>n^{1/3}}-\expec[W^2\indic{W>n^{1/3}}]\Big|>\vep n^{-1/3}\Big)\\
    &\qquad\leq n^{1/3}\expec[W^2\indic{W>n^{1/3}}]
    \leq\expec[W^ 3\indic{W>n^{1/3}}]=o(1),\nn
    }
since $\expec[W^3]<\infty$. Thus,
    \eqn{
    \frac{n^{1/3}}{n}\sum_{i=1}^n \Big|W_i^2\indic{W_i>n^{1/3}}-\expec[W^2\indic{W>n^{1/3}}]\Big|
    \convp 0.
    }
For the first sum, we use the Chebycheff inequality to obtain
    \eqan{
    \prob\Big(\Big|\frac{1}{n}\sum_{i=1}^n W_i^2\indic{W_i\leq n^{1/3}}-\expec[W^2]\Big|>\vep n^{-1/3}\Big)
    &=\vep^{-2} n^{2/3}{\rm Var}\Big(\frac{1}{n}\sum_{i=1}^n W_i^2\indic{W_i\leq n^{1/3}}\Big)\nn\\
    &\leq \vep^{-2} n^{-1/3} \expec\Big[W^4 \indic{W\leq n^{1/3}}\Big]=o(1),
    }
since, when $\expec[W^3]<\infty$, we have that $\expec[W^4\indic{W\leq x}]=o(x)$.
Thus, also,
    \eqn{
    \frac{n^{1/3}}{n}\sum_{i=1}^n \Big|W_i^2\indic{W_i\leq n^{1/3}}-\expec[W^2\indic{W\leq n^{1/3}}]\Big|
    \convp 0.
    }
This proves that conditions (b)-(c) hold in probability.

\subsection{Verification of conditions for weights as in \refeq{choicewi}}
Here we check conditions (b) and (c) for the case that
$\bfw=(w_1, \ldots, w_n)$ where $w_i$ is chosen as in
\refeq{choicewi}. We shall frequently make use of the fact that
\refeq{[1-F]bd} implies that $1-F(x)=o(x^{-3})$ as $x\rightarrow \infty$,
which, in turn implies that (see e.g., \cite[(B.9)]{EskHofHoo06}), as $u\downarrow 0$,
    \eqn{
    \lbeq{[1-F]inv-bd}
    [1-F]^{-1}(u)=o(u^{-1/3}).
    }
To verify condition (b), we note that by
\cite[(4.2)]{Hofs09a}, $w_{V_n}$ has distribution function
    \eqn{
    \lbeq{Fn-def}
    F_n(x)=\frac 1n\big(\big\lfloor n F(x)\big\rfloor+1\big)\wedge 1.
    }
This converges to $F(x)$ for every $x\geq 0$, which proves that condition (b)
holds. To verify condition (c), we note that since $i\mapsto [1-F]^{-1}(i/n)$
is monotonically decreasing, for any $s>0$, we have
    \eqn{
    \expec[W^s]-\int_0^{1/n} [1-F^{-1}(u)]^s du \leq \frac{1}{n}\sum_{i=1}^n w_i^s\leq \expec[W^s].
    }
Now, by \refeq{[1-F]inv-bd}, we have that, for $s=1, 2, 3$,
    \eqn{
    \int_0^{1/n} [1-F^{-1}(u)]^s du=o(n^{s/3-1}),
    }
which proves all necessary bounds
for condition (c) at once.

\paragraph{Acknowledgements.}
We thank Tatyana Turova for lively and open discussions, and for inspiring us to 
present our results for a more general setting. The work of SSB was supported by 
PIMS and NSERC, Canada. The work of RvdH and JvL work was supported in part by 
the Netherlands Organisation for Scientific Research (NWO).

\def\cprime{$'$}



\begin{thebibliography}{10}

\bibitem{Aldo97}
D.~Aldous.
\newblock Brownian excursions, critical random graphs and the multiplicative
  coalescent.
\newblock {\em Ann. Probab.}, {\bf 25}(2):812--854, (1997).

\bibitem{aldous-limic}
D.~Aldous and V.~Limic.
\newblock {The entrance boundary of the multiplicative coalescent}.
\newblock {\em Electron. J. Probab}, 3:1--59, 1998.

\bibitem{BhaHofLee09b}
S.~Bhamidi, R.~van~der Hofstad, and J.S.H.~van Leeuwaarden.
\newblock Novel scaling limits for critical inhomogeneous random graphs.
\newblock Preprint (2009).

\bibitem{Boll01}
B.~Bollob{\'a}s.
\newblock {\em Random graphs}, volume~{\bf 73} of {\em Cambridge Studies in
  Advanced Mathematics}.
\newblock Cambridge University Press, Cambridge, second edition, (2001).

\bibitem{BolJanRio07}
B.~Bollob{\'a}s, S.~Janson, and O.~Riordan.
\newblock The phase transition in inhomogeneous random graphs.
\newblock {\em Random Structures Algorithms}, {\bf 31}(1):3--122, (2007).

\bibitem{BriDeiMar-Lof05}
T.~Britton, M.~Deijfen, and A.~Martin-L{\"o}f.
\newblock Generating simple random graphs with prescribed degree distribution.
\newblock {\em J. Stat. Phys.}, {\bf 124}(6):1377--1397, (2006).

\bibitem{ChuLu02a}
F.~Chung and L.~Lu.
\newblock The average distances in random graphs with given expected degrees.
\newblock {\em Proc. Natl. Acad. Sci. USA}, {\bf 99}(25):15879--15882
  (electronic), (2002).

\bibitem{ChuLu02b}
F.~Chung and L.~Lu.
\newblock Connected components in random graphs with given expected degree
  sequences.
\newblock {\em Ann. Comb.}, {\bf 6}(2):125--145, (2002).

\bibitem{ChuLu03}
F.~Chung and L.~Lu.
\newblock The average distance in a random graph with given expected degrees.
\newblock {\em Internet Math.}, {\bf 1}(1):91--113, (2003).

\bibitem{ChuLu06c}
F.~Chung and L.~Lu.
\newblock {\em Complex graphs and networks}, volume~{\bf 107} of {\em CBMS
  Regional Conference Series in Mathematics}.
\newblock Published for the Conference Board of the Mathematical Sciences,
  Washington, DC, (2006).

\bibitem{ChuLu06}
F.~Chung and L.~Lu.
\newblock The volume of the giant component of a random graph with given
  expected degrees.
\newblock {\em SIAM J. Discrete Math.}, {\bf 20}:395--411, (2006).

\bibitem{EskHofHoo06}
H.~van~den Esker, R.~van~der Hofstad, and G.~Hooghiemstra.
\newblock Universality for the distance in finite variance random graphs.
\newblock {\em J. Stat. Phys.}, {\bf 133}(1):169--202, (2008).

\bibitem{EthKur86}
S.N. Ethier and T.G. Kurtz.
\newblock {\em Markov processes}.
\newblock Wiley Series in Probability and Mathematical Statistics: Probability
  and Mathematical Statistics. John Wiley \& Sons Inc., New York, (1986).
\newblock Characterization and convergence.

\bibitem{Hofs09a}
R.~van~der Hofstad.
\newblock Critical behavior in inhomogeneous random graphs.
\newblock Preprint (2009).

\bibitem{HofKagMul09}
R.~van~der Hofstad, W.~Kager, and T.~M{\"u}ller.
\newblock A local limit theorem for the critical random graph.
\newblock {\em Electron. Commun. Probab.}, {\bf 14}:122--131, (2009).

\bibitem{Jans08a}
S.~Janson.
\newblock Asymptotic equivalence and contiguity of some random graphs.
\newblock Preprint (2008).

\bibitem{JanLucRuc00}
S.~Janson, T.~{\L}uczak, and A.~Rucinski.
\newblock {\em Random graphs}.
\newblock Wiley-Interscience Series in Discrete Mathematics and Optimization.
  Wiley-Interscience, New York, (2000).

\bibitem{Mart98}
A.~Martin-L{\"o}f.
\newblock The final size of a nearly critical epidemic, and the first passage
  time of a {W}iener process to a parabolic barrier.
\newblock {\em J. Appl. Probab.}, {\bf 35}(3):671--682, (1998).

\bibitem{NacPer07a}
A.~Nachmias and Y.~Peres.
\newblock Critical percolation on random regular graphs.
\newblock Preprint (2007).

\bibitem{NorRei06}
I.~Norros and H.~Reittu.
\newblock On a conditionally {P}oissonian graph process.
\newblock {\em Adv. in Appl. Probab.}, {\bf 38}(1):59--75, (2006).

\bibitem{Pitt01}
B.~Pittel.
\newblock On the largest component of the random graph at a nearcritical stage.
\newblock {\em J. Combin. Theory Ser. B}, {\bf 82}(2):237--269, (2001).

\bibitem{Turo09}
T.S. Turova.
\newblock Diffusion approximation for the components in critical inhomogeneous
  random graphs of rank 1.
\newblock (2009).

\end{thebibliography}
\end{document}